\documentclass[12pt]{article}

\usepackage{amsfonts}
\usepackage{graphicx} % for \includegraphics
\usepackage{subcaption} % for subfigure environment
\usepackage{amsmath} % for math environment
\usepackage{listings} % for code environment
\usepackage{xcolor} % for different text colors
\usepackage{hyperref} % for hyperlinks
\usepackage{fancyhdr} % for header and footer
\usepackage{amssymb} % for mathbb and others
\usepackage{amsthm} % proofs, theorems, lemmas, remarks...
\usepackage{mathtools}
\usepackage{authblk}
\usepackage{float} % force figure placement

\newtheorem{theorem}{Theorem}

\newtheorem{definition}[theorem]{Definition}
\newtheorem{example}[theorem]{Example}

\newtheorem{lemma}[theorem]{Lemma}

\hypersetup{
colorlinks=true,
linkcolor=blue, filecolor=magenta, urlcolor=cyan, }

\title{Periodic second-order systems and coupled forced Van der Pol oscillators}
\author[1,2]{Feliz Minh\'{o}s}
\author[1*]{Sara Perestrelo}
\affil[1]{Department of Mathematics, School of Sciences and Technology, University of \'{E}vora, Rua Rom\~{a}o
Ramalho, 59, 7000-671 \'{E}vora, Portugal}
\affil[2]{Research Center in Mathematics and Applications (CIMA), Institute for Advanced Studies and Research (IIFA), University of \'{E}vora, Portugal}
\affil[*]{Correspondence: belperes@uevora.pt}
\date{}               %% if you don't need date to appear
\providecommand{\keywords}[1]
{\small	\textbf{Keywords:} #1}

\providecommand{\MSC}[1]
{\small	\textbf{MSC:} #1}

\begin{document}
\maketitle

\begin{abstract}
We present an existence and localization result for periodic solutions of second-order non-linear coupled planar systems, without requiring periodicity for the non-linearities. The arguments for the existence tool are based on a variation of the Nagumo condition and the Topological Degree Theory. The localization tool is based on a technique of orderless upper and lower solutions, that involves functions with translations. We apply our result to a system of two coupled Van der Pol oscillators with a forcing component.
\end{abstract}


\textbf{2020} \MSC{34A34, 34A12, 34C25, 34C15.}

\bigskip

\section{Introduction}
\label{sec: introduction}

\quad We study the following second-order non-linear coupled system,
\begin{equation}
\left\{ 
\begin{array}{c}
z^{\prime \prime }(t)=f\left(t, z(t), w(t), z'(t), w'(t)\right) \\ 
w^{\prime \prime }(t)=g\left(t, z(t), w(t), z'(t), w'(t)\right)%
\end{array}%
\right. ,  \label{eq: OP}
\end{equation}%
for $t\in \left[ 0,T\right]$, $T>0$, with continuous functions $f,g:\left[ 0,T\right] \times 
\mathbb{R}^{4}\rightarrow \mathbb{R}$, and the periodic boundary conditions 
\begin{equation}
\begin{array}{ccc}
&z(0) = z(T),\qquad  &z'(0)=z'(T),\\
&w(0) = w(T),\qquad &w'(0)=w'(T).
\end{array}
\label{eq: PBC}
\end{equation}

Coupled ordinary differential systems arise in many domains of life, namely Chemistry \cite{varadharajan2011analytical}, Biology \cite{julius2008stochastic}, Physiology \cite{mahata2017application}, Mechanics \cite{oliveira2022stability}, Economy \cite{ikeda2012coupled}, Electronics \cite{orihashi2005one}, among others.

Finding periodic solutions in second-order differential equations or systems can be delicate, specially when the non-linearities are non-periodic. The search for periodic solutions in such problems has had an increasing interest, as it can be seen by several papers from the last decades, for example, \cite{lu2019periodic, wang2012ambrosetti, fonda2021periodic}. The most common argument for obtaining periodic solutions is based on the assumptions of periodicity in the non-linearities. As an example, we refer to the work \cite{wang2012ambrosetti}, where the authors study the problem
\begin{equation}
\left\{
\begin{split}
& u_{tt} - u_{xx} + c_1 u_t + \Xi_1(t, x, u, v) = s_1,\\
& v_{tt} - v_{xx} + c_2 v_t + \Xi _2(t, x, u, v) = s_2,\\
& u(t + 2 \pi, x) = u(t, x + 2 \pi) = u(t, x),\\
& v(t + 2 \pi, x) = v(t, x + 2 \pi) = v(t, x),
\end{split}
\right.
\label{eq: intro - ex1}
\end{equation}
with constants $c_i>0$, parameters $s_i$, continuous functions $\Xi_i(t, x, u, v): \mathbb{R}^4 \rightarrow \mathbb{R}$, $2\pi$-periodic in $t$ and $x$, obtaining an Ambrosetti-Prodi alternative for periodic solutions.

In \cite{lu2019periodic}, the following second-order differential equation is studied,
\begin{equation}
\chi''(t) + \Psi(\chi(t))\,\chi'(t) + \phi(t) \chi^m(t) - \frac{a(t)}{\chi^{\mu}(t)} + \frac{b(t)}{\chi^{y}(t)},
\label{eq: intro - ex2}
\end{equation}
where $\Psi \in C((0,+\infty), \mathbb{R})$, $\phi$, $a$ and $b$ are $T$-periodic, and in $L([0,T], \mathbb{R})$, $m$, $\mu$ and $y$ are constants with $m \geq 0$ and $\mu \geq y > 0$. Here, asymptotic arguments are a key argument to obtain periodic solutions.

The above problems, (\ref{eq: intro - ex1}) and (\ref{eq: intro - ex2}), require periodicity for the non-linearities. In \cite{fonda2021periodic}, the authors study the second-order differential problem,
\begin{equation}
\left\{
\begin{split}
& \ddot{x}(t) = \Theta(t, x),\\
& x(0) = x(T), \\
& \dot{x}(0) = \dot{x}(T),
\end{split}
\right.
\label{eq: intro - ex3}
\end{equation}
without requiring periodicity for the non-linearities. However, the function $\Theta$ does not depend on the first derivative, thus making problem (\ref{eq: intro - ex3}) a particular case of problem (\ref{eq: OP}). Moreover, the obtained results require an order relation (for the well and inverse order cases) between lower and upper solutions.

Other attempts to find periodic solutions in ordinary differential problems may not present periodicity requirements but, typically, only approximate solutions are found, see \cite{yildirim2011higher, ju2014global, mickens2005generalized}.

To overcome such constraints, in this work we present an original methodology for finding periodic solutions in generalized second-order coupled systems.

Motivated by previous works existent in the literature, (see, for instance, \cite{minhos2019solvability, minhos2020periodic, minhos2023periodic, minhos2023coupled}), it is the first time, to the best of our knowledge, where sufficient conditions are given to prove the existence of periodic solutions for coupled systems, without assuming periodicity in the non-linearities nor any order requirement between lower and upper solutions. Therefore, the set of admissible functions for possible lower and upper solutions has a wider range.

Within this line of research, the present work is somehow a continuation of \cite{minhs2023first, minhos2023first}, which aim to provide a complete study on the existence and localization of periodic solutions in first-order generalized coupled systems, with and without impulses, using this technique of orderless lower and upper solutions, where the lack of order is overcome with adequate translations.

The strategy employed consists of proving the existence of at least a periodic solution $(z(t), w(t))$ applying the Topological Degree Theory and to localize the existing solution in a strip bounded by lower and upper functions.

This paper is organized as follows. In Section \ref{sec: definitions} all the required definitions and lemmas are presented. In Section \ref{sec: result} the existence and localization result is formulated as the main theorem, together with the respective proof in four claims. A numerical example shows the applicability of the main theorem. In Section \ref{sec: application} we apply our methodology to the problem of two coupled Van der Pol oscillators with forcing terms. As all the requisites of the main theorem in Section \ref{sec: result} are verified, we suggest possible localizing functions for the existing periodic solutions of this problem. We draw the conclusions of our study in Section \ref{sec: conclusions}.

\bigskip

\section{Definitions}
\label{sec: definitions}

\quad Consider the Banach space $X:=C^1[0,T]$, equipped with the norm
\begin{equation*}
\Vert u \Vert_X := \max\{ \Vert u \Vert, \Vert u' \Vert \}, \qquad \Vert u \Vert := \max_{t \in [0,T]} \{ \vert u(t) \vert \},
\end{equation*}
and the vectorial space $X^2 :=(C^1[0,T])^2$, with the norm
\begin{equation*}
\Vert (u,v) \Vert_{X^2} = \max \{ \Vert u \Vert_X, \Vert v \Vert_X \}.
\end{equation*}

\bigskip

In second-order differential equations it is important to have a control on the first derivatives, which is achieved using a Nagumo-type condition:

\begin{definition}\label{def: nagumo}
Consider $C^1$ continuous functions $\gamma_i, \Gamma_i: [0,T]\rightarrow\mathbb{R}, i=1,2$, and the set
\begin{equation*}
S = \{ (t, z_0, w_0, z_1, w_1) \in [0,T]\times \mathbb{R}^4: \gamma_1(t)\leq z_0 \leq \Gamma_1(t), \gamma_2(t)\leq w_0 \leq \Gamma_2(t) \}.
\end{equation*}
The continuous functions $f,g: [0,T]\times \mathbb{R}^4\rightarrow\mathbb{R}$  satisfy a Nagumo-type condition relative to the intervals $[\gamma_1(t), \Gamma_1(t)]$ and $[\gamma_2(t), \Gamma_2(t)]$, for all $t \in [0,T]$, if there exist continuous functions $\varphi, \psi: [0, +\infty[ \rightarrow ]0, +\infty[$ verifying
\begin{equation}
\int_{0}^{+\infty} \frac{ds}{\varphi(\vert s \vert)} = + \infty,
\qquad
\int_{0}^{+\infty} \frac{ds}{\psi(\vert s \vert)} = + \infty,
\label{eq: phi,psi}
\end{equation}
such that
\begin{equation}
\begin{split}
\vert f(t, z_0, w_0, z_1, w_1) \vert  \leq \varphi(\vert z_1 \vert),
\quad \forall (t, z_0, w_0, z_1, w_1) \in S, \\
\vert g(t, z_0, w_0, z_1, w_1) \vert  \leq \psi(\vert w_1 \vert),
\quad \forall (t, z_0, w_0, z_1, w_1) \in S.
\end{split}
\label{eq: f,g < phi, psi}
\end{equation}
\end{definition}

The \textit{a priori} estimates for the first derivatives are given by the following lemma.

\begin{lemma}\label{lemma: estimates}
Suppose that the continuous functions $f,g:[0,T]\times \mathbb{R}^4 \rightarrow \mathbb{R}$ satisfy a Nagumo-type condition relative to the intervals $[\gamma_1(t), \Gamma_1(t)]$ and $[\gamma_2(t), \Gamma_2(t)]$, for all $t \in [0,T]$.

Then, for every solution $(z(t),w(t)) \in X^2$ of (\ref{eq: OP}), (\ref{eq: PBC}) verifying
\begin{equation}
\gamma_1(t) \leq z(t) \leq \Gamma_1(t),
\quad
\textrm{and}
\quad
\gamma_2(t) \leq w(t) \leq \Gamma_2(t),
\quad
\forall t \in [0,T],
\label{eq: lemma, z,w localization}
\end{equation}
there are $N_1, N_2 > 0 $ such that
\begin{equation}
\Vert z' \Vert \leq N_1, \qquad\Vert w' \Vert \leq N_2.
\label{eq: N1, N2}
\end{equation}
\end{lemma}

\begin{proof}
Let $(z(t), w(t))$ be a solution of (\ref{eq: OP}), (\ref{eq: PBC}) verifying (\ref{eq: lemma, z,w localization}).

By (\ref{eq: PBC}) and by Rolle's Theorem, there exists $t_0\in[0,T]$ such that $z'(t_0)=0$.

Consider $N_i>0, i=1,2$ such that
\begin{equation}
\int_0^{N_1} \frac{ds}{\varphi(\vert s \vert)} > T,
\qquad
\int_0^{N_2} \frac{ds}{\psi(\vert s \vert)} > T,
\label{eq: Ni}
\end{equation}
and assume, without loss of generality, there exist $t_1, t_2 \in [0,T]$ such that $z'(t_1) \leq 0$ and $z'(t_2)>0$.

By continuity of $z'(t)$, there exists $t_3 \in [t_1, t_2]$ such that $z'(t_3)=0$. Using a convenient change of variables, and by (\ref{eq: OP}) and (\ref{eq: f,g < phi, psi}),
\begin{equation*}
\begin{split}
\int^{z'(t_2)}_{z'(t_3)} \frac{ds}{\varphi(\vert s \vert)}
&=
\int^{t_2}_{t_3} \frac{z''(t)}{\varphi(\vert z'(t) \vert)} dt
\leq
\int_0^T \frac{\vert z''(t) \vert}{\varphi(\vert z'(t) \vert)} dt \\[5pt]
&=
\int_0^T \frac{\vert f(t,z(t),w(t),z'(t),w'(t)) \vert}{\varphi(\vert z'(t) \vert)} dt
\leq
\int_0^T \frac{\varphi(\vert z'(t) \vert)}{\varphi(\vert z'(t) \vert)} dt = T.
\end{split}
\end{equation*}
By (\ref{eq: Ni}),
\begin{equation}
\int^{z'(t_2)}_{z'(t_3)} \frac{ds}{\varphi(\vert s \vert)}
=
\int^{N_1}_{0} \frac{ds}{\varphi(\vert s \vert)}
\leq
T
<
\int^{N_1}_{0} \frac{ds}{\varphi(\vert s \vert)},
\end{equation}
and, therefore, $z'(t_2)<N_1$. As $t_2$ is chosen arbitrarily, then $z'(t)<N_1$, for all $t \in [0,T]$.

The case where $z'(t_1)>0$ and $z'(t_2)\leq 0$ follows similar arguments. Therefore, $\Vert z' \Vert \leq N_1$, for all $t \in [0,T]$. Likewise, $\Vert w'\Vert \leq N_2$, for all $t \in [0,T]$.

Remark that the constant $N_1$ depends only on $\gamma_i$, $\Gamma_i$, $\varphi$ and $T$. In the same way, $N_2$ depends only on $\gamma_i$, $\Gamma_i$, $\psi$ and $T$.

\end{proof}

We present below the definition of upper and lower solutions.

\begin{definition}\label{def: lower upper}
The pair of real functions $\left(\alpha_{1},\alpha_{2}\right) \in (C^1[0,T])^2$ is a lower solution of the periodic problem (\ref{eq: OP}), (\ref{eq: PBC}) if
\begin{equation*}
\begin{split}
&\alpha_1''(t) \geq f(t, \alpha_1^0(t), \alpha_2^0(t), \alpha_1'(t), w_1), \forall w_1 \in \mathbb{R},\\
&\alpha_2''(t) \geq g(t, \alpha_1^0(t), \alpha_2^0(t), z_1, \alpha_2'(t)), \forall z_1 \in \mathbb{R},
\end{split}
\end{equation*}
with
\begin{equation*}
\alpha_i^0(t):= \alpha_i(t)-\Vert \alpha_i \Vert,
\end{equation*}
and
\begin{equation*}
\alpha_i(0) = \alpha_i(T), \qquad \alpha'_i(0) \geq \alpha'_i(T).
\end{equation*}

The pair of real functions $\left(\beta_{1},\beta_{2}\right) \in (C^1[0,T])^2$ is an upper solution of the periodic problem (\ref{eq: OP}), (\ref{eq: PBC}) if
\begin{equation}
\begin{split}
&\beta_1''(t) \leq f(t, \beta_1^0(t), \beta_2^0(t), \beta_1'(t), w_1), \forall w_1 \in \mathbb{R},\\
&\beta_2''(t) \leq g(t, \beta_1^0(t), \beta_2^0(t), z_1, \beta_2'(t)), \forall z_1 \in \mathbb{R},
\end{split}
\label{eq: beta'' < f lower upper}
\end{equation}
with
\begin{equation}
\beta_i^0(t):= \beta_i(t)+\Vert \beta_i \Vert,
\label{eq: beta0}
\end{equation}
and
\begin{equation}
\beta_i(0) = \beta_i(T), \qquad \beta_i'(0) \leq \beta_i'(T).
\label{eq: beta - monotonicity}
\end{equation}
\end{definition}

\bigskip

\section{Main Result}
\label{sec: result}

\quad The following theorem is an existence and localization result using orderless lower and upper solutions.

\begin{theorem}\label{thm: 1}
Let $f,g:[0,T] \times \mathbb{R}^4 \rightarrow \mathbb{R}$ be continuous functions. Suppose $(\alpha_1,\alpha_2)$, $(\beta_1,\beta_2)$ are lower and upper solutions of problem (\ref{eq: OP}), (\ref{eq: PBC}), as in Definition \ref{def: lower upper}.

Suppose $f,g$ satisfy a Nagumo-type condition, according to Definition \ref{def: nagumo} relative to the intervals $[\alpha_1^0(t), \beta_1^0(t)]$ and $[\alpha_2^0(t), \beta_2^0(t)]$, for all $t \in [0,T]$, with
\begin{equation}\label{eq: f non-increasing}
f(t, z_0,w_0,z_1,w_1) \textrm{ non-increasing in }w_0,
\textrm{ for } t \in [0,T], 
z_0 \in \mathbb{R}\textrm{ fixed},
\end{equation}
\begin{equation}
\min\left\{ \min_{t\in [0,T]} \alpha_1'(t), \min_{t\in [0,T]} \beta_1'(t) \right\}
\leq
z_1
\leq
\max\left\{ \max_{t\in [0,T]} \alpha_1'(t), \min_{t\in [0,T]} \beta_1'(t) \right\},
\end{equation}
and
\begin{equation*}
g(t, z_0,w_0,z_1,w_1) \textrm{ non-increasing in }z_0,
\textrm{ for } t \in [0,T], 
w_0 \in \mathbb{R}\textrm{ fixed},
\end{equation*}
\begin{equation*}
\min\left\{ \min_{t\in [0,T]} \alpha_2'(t), \min_{t\in [0,T]} \beta_2'(t) \right\}
\leq
w_1
\leq
\max\left\{ \max_{t\in [0,T]} \alpha_2'(t), \min_{t\in [0,T]} \beta_2'(t) \right\}.
\end{equation*}

Then, there exists at least one pair $(z(t), w(t)) \in (C^2[0,T])^2$, solution of problem (\ref{eq: OP}), (\ref{eq: PBC}) and, moreover,
\begin{equation*}
\alpha_1^0(t) \leq z(t) \leq \beta_1^0(t),
\qquad
\alpha_1^0(t) \leq w(t) \leq \beta_1^0(t),
\qquad
\forall t \in [0,T].
\end{equation*}
\end{theorem}

\bigskip

\begin{proof}
Define the truncated functions $\delta_i:[0,T]\times\mathbb{R}\rightarrow\mathbb{R}$,
\begin{equation}
\begin{array}{ccc}
&\delta_1(t,z)=
\left\{
\begin{split}
&\beta_1^0(t), \quad  z>\beta_1^0(t) \\
&z, \quad  \alpha_1^0(t)\leq z \leq \beta_1^0(t) \\
&\alpha_1^0(t), \quad  z<\alpha_1^0(t) \\
\end{split}
\right.,
\quad
&\delta_2(t,w)=
\left\{
\begin{split}
&\beta_2^0(t), \quad  w>\beta_2^0(t) \\
&w, \quad  \alpha_2^0(t)\leq w \leq \beta_2^0(t) \\
&\alpha_2^0(t), \quad  w<\alpha_2^0(t) \\
\end{split}
\right..
\end{array}
\end{equation}
For $\lambda, \mu \in [0,1]$, consider the truncated, perturbed and homotopic auxiliary problem
\begin{equation}
\left\{
\begin{split}
&z''(t)-z(t)=
\lambda \, f(t, \delta_1(t,z(t)), \delta_2(t,w(t)), z'(t), w'(t))-\lambda \, \delta_1(t, z(t)) \\
&w''(t)-w(t)=
\mu \, g(t, \delta_1(t,z(t)), \delta_2(t,w(t)), z'(t), w'(t))-\mu \, \delta_2(t, w(t)) 
\end{split}
\right.,
\label{eq: AP}
\end{equation}
for all $t \in [0,T]$, with the boundary conditions (\ref{eq: PBC}).

Consider $r_i>0, i=1,2$, such that, for all $\lambda, \mu \in [0,1]$, every $t \in [0,T]$, and every $N_1^*,N_2^* >0$ given by (\ref{eq: N1, N2}),
\begin{equation}
-r_i < \alpha_i^0(t)\leq\beta_i^0(t)<r_i,
\label{eq: ri}
\end{equation}
with
\begin{equation}\label{eq: condicoes alphas betas ri}
\begin{split}
&\lambda \, f(t, \beta_1^0(t), \beta_2^0(t), 0, w'(t))- \beta_1^0(t)+r_1 >0,
\quad \textrm{ for } \Vert w' \Vert \leq N_2^*,\\
&\mu \, g(t, \beta_1^0(t), \beta_2^0(t), z'(t), 0)- \beta_2^0(t)+r_2 >0,
\quad \textrm{ for } \Vert z' \Vert \leq N_1^*,\\
&\lambda \, f(t, \alpha_1^0(t), \alpha_2^0(t), 0, w'(t))- \alpha_1^0(t)-r_1 <0,\quad \textrm{ for } \Vert w' \Vert \leq N_2^*,\\
&\mu \, g(t, \alpha_1^0(t), \alpha_2^0(t), z'(t), 0)- \alpha_2^0(t)-r_2 <0,
\quad \textrm{ for } \Vert z' \Vert \leq N_1^*.
\end{split}
\end{equation}

\bigskip

\textbf{Claim 1.} \textit{Every solution of (\ref{eq: AP}), (\ref{eq: PBC}) verifies $\vert z(t) \vert < r_1$ and $\vert w(t) \vert < r_2$, for all $t \in [0,T]$, independently of $\lambda, \mu \in [0,1]$.}

Assume, by contradiction, that there exist $\lambda \in [0,1]$, a pair $(z(t), w(t))$, solution of problem (\ref{eq: AP}), (\ref{eq: PBC}), and $t \in [0,T]$ such that $\vert z(t) \vert \geq r_1$. If $z(t) \geq r_1$, define
\begin{equation}
\max_{t\in[0,T]} z(t):=z(t_0).
\end{equation}

If $t\in]0,T[$ and $\lambda \in]0,1]$, then $z'(t)=0$ and $z''(t)\leq 0$. By (\ref{eq: ri}), (\ref{eq: f non-increasing}) and (\ref{eq: condicoes alphas betas ri}), the following contradiction holds:
\begin{equation*}
\begin{split}
0\geq z''(t_0)
&=
\lambda \, f(t_0, \delta_1(t_0,z(t_0)), \delta_2(t_0,w(t_0)), z'(t_0), w'(t_0))-\lambda \, \delta_1(t_0, z(t_0)) + z(t_0)\\
&=
\lambda \, f(t_0, \beta_1^0(t_0), \delta_2(t_0,w(t_0)), 0, w'(t_0))-\lambda \, \beta_1^0(t_0) + z(t_0)\\
& \geq
\lambda \, f(t_0, \beta_1^0(t_0), \delta_2(t_0,w(t_0)), 0, w'(t_0))-\beta_1^0(t_0) + z(t_0)\\
& \geq
\lambda \, f(t_0, \beta_1^0(t_0), \delta_2(t_0,w(t_0)), 0, w'(t_0))-\beta_1^0(t_0) + r_1
>0
\end{split}
\end{equation*}

If $t_0=0$ or $t_0=T$, then, by (\ref{eq: PBC}),
\begin{equation*}
0 \geq z'(0)=z'(T) \geq 0.
\end{equation*}
So, $z'(0)=z'(T)=0$, $z''(0)\leq 0$ and $z''(T)\leq 0$. Therefore, the arguments follow the previous case.

If $\lambda=0$, we obtain the contradiction
\begin{equation*}
0 \geq z''(t_0)=z(t_0)\geq r_1>0.
\end{equation*}

Then, $z(t)<r_1, \forall t \in [0,T]$, regardless of $\lambda$.

The same arguments can be made to prove that $z(t)>-r_1$. Therefore, $\vert z(t) \vert < r_1$, for $t \in [0,T]$, independently of $\lambda$.

Likewise, $\vert w(t) \vert < r_2$, for $t \in [0,T]$, independently of $\mu$.

\bigskip

\textbf{Claim 2.} \textit{For every solution $(z(t),w(t))$ of (\ref{eq: AP}), (\ref{eq: PBC}), there are $N_1^*, N_2^*>0$ such that $\vert z'(t) \vert < N_1^*$ and $\vert w'(t) \vert < N_2^*$, $\forall t \in [0,T]$, independently of $\lambda, \mu \in [0,1]$.}

Define the continuous functions
\begin{equation*}
\begin{split}
&F_{\lambda}(t,z_0,w_0,z_1,w_1)):=
\lambda \, f(t,\delta_1(t,z_0),\delta_2(t,w_0),z_1,w_1) - \lambda \, \delta_1(t,z_0) + z_0,\\
&G_{\mu}(t,z_0,w_0,z_1,w_1):=
\mu \, f(t,\delta_1(t,z_0),\delta_2(t,w_0),z_1,w_1) - \mu \, \delta_2(t,w_0) + w_0,
\end{split}
\end{equation*}
and, as $f,g:[0,T]\times\mathbb{R}^4\rightarrow\mathbb{R}$ satisfy a Nagumo-type condition relative to the intervals $[\alpha_1^0(t),\beta_1^0(t)]$ and $[\alpha_2^0(t),\beta_2^0(t)]$, then,
\begin{equation*}
\begin{split}
&\vert F_{\lambda}(t,z(t),w(t),z'(t),w'(t)) \vert \\
&\qquad\leq
\vert f(t,\delta_1(t,z(t)),\delta_2(t,w(t)),z'(t),w'(t)) \vert +\vert \delta_1(t,z(t)) \vert + \vert z(t) \vert\\
&\qquad < \varphi(\vert z' \vert) + 2r_1,\\
&\vert G_{\mu}(t,z(t),w(t),z'(t),w'(t)) \vert \\
&\qquad\leq
\vert f(t,\delta_1(t,z(t)),\delta_2(t,w(t)),z'(t),w'(t)) \vert +\vert \delta_2(t,w(t)) \vert + \vert w(t) \vert\\
&\qquad < \psi(\vert w' \vert) + 2r_2.
\end{split}
\end{equation*}

Therefore, for continuous positive functions $\varphi^*,\psi^*:[0,+\infty[\rightarrow]0,+\infty[$ given by
\begin{equation*}
\varphi^*(\vert z' \vert):=\varphi(\vert z' \vert) + 2r_1, \qquad \psi^*(\vert w' \vert):=\psi(\vert w' \vert) + 2r_2,
\end{equation*}
$F_{\lambda}$ and $G_{\mu}$ satisfy a Nagumo-type condition in 
\begin{equation*}
E=\{ (t,z,w,z',w') \in [0,T]\times\mathbb{R}^4:\vert z \vert <r_1, \vert w \vert <r_2 \},
\end{equation*}
as, by (\ref{eq: phi,psi}), we have
\begin{equation}
\begin{split}
&\int_0^{+\infty}\frac{ds}{\varphi^*(\vert z' \vert)} = 
\int_0^{+\infty}\frac{ds}{\varphi(\vert z' \vert)+2r_1}=+\infty, \\
&\int_0^{+\infty}\frac{ds}{\psi^*(\vert w' \vert)} = 
\int_0^{+\infty}\frac{ds}{\psi(\vert w' \vert)+2r_2}=+\infty.
\end{split}
\end{equation}

Therefore, by Lemma \ref{lemma: estimates}, there are $N_1^*,N_2^*>0$ such that $\vert z'(t) \vert < N_1^*$ and $\vert w'(t) \vert < N_2^*$, $\forall t \in [0,T]$, independently of $\lambda, \mu \in [0,1]$.

\bigskip

\textbf{Claim 3.} \textit{Problem (\ref{eq: AP}), (\ref{eq: PBC}), for $\lambda=\mu=1$, has at least one solution $(z(t), w(t))$.}

Define the operators
\begin{equation*}
\mathcal{L}:(C^2[0,T])^2 \rightarrow (C^1[0,T])^2 \times \mathbb{R}^4, 
\end{equation*}
given by 
\begin{equation*}
\mathcal{L} (z, w) := (z''(t)-z(t), w''(t)-w(t), z(0), w(0), z'(0), w'(0)),
\end{equation*}
and
\begin{equation*}
\mathcal{N}_{(\lambda, \mu)}: (C^2[0,T])^2 \rightarrow (C^1[0,T])^2 \times \mathbb{R}^4, 
\end{equation*}
given by
\begin{equation*}
\mathcal{N}_{(\lambda,\mu)} (z, w) := (Z_{\lambda}(t), W_{\mu}(t), z(T), w(T), z'(T), w'(T)),
\end{equation*}
where
\begin{equation*}
\begin{split}
& Z_{\lambda}(t):=
\lambda \, f(t,\delta_1(t,z(t)),\delta_2(t,w(t)),z'(t),w'(t)) - \lambda \, \delta_1(t,z(t)),\\
& W_{\mu}(t):=
\mu \, f(t,\delta_1(t,z(t)),\delta_2(t,w(t)),z'(t),w'(t)) - \mu \, \delta_2(t,w(t)).
\end{split}
\end{equation*}

We define the continuous operator
\begin{equation*}
\mathcal{T}: (C^2[0,T])^2 \rightarrow (C^2[0,T])^2, 
\end{equation*}
given by
\begin{equation*}
\mathcal{T}_{(\lambda,\mu)} (z, w):= \mathcal{L}^{-1}\mathcal{N}_{(\lambda, \mu)}(z, w).
\end{equation*}
For $M:=\max\{ r_1, r_2, N_1^*, N_2^* \}$, consider the set
\begin{equation}
\Omega = \{ (z,w) \in (C^2([0,T]))^2 : \Vert (z,w) \Vert_{X^2} < M \}.
\end{equation}

By Claims 1 and 2, for all $\lambda, \mu \in [0,1]$, the degree $d(I-\mathcal{T}_{(0,0)}, \Omega, 0)$ is well defined and, by homotopy invariance,
\begin{equation}
d(I-\mathcal{T}_{(0,0)}, \Omega, 0)
=
d(I-\mathcal{T}_{(1,1)}, \Omega, 0).
\label{eq: homotopy invariance}
\end{equation}
As the equation $(z,w)=\mathcal{T}_{(0,0)}(z,w)$, that is, the homogeneous system
\begin{equation*}
\left\{
\begin{split}
& z''(t)-z(t)=0\\
& w''(t)-w(t)=0,
\end{split}
\right.
\end{equation*}
with the periodic conditions (\ref{eq: PBC}), admits only the null solution, then $\mathcal{L}^{-1}$ is well-defined and, by degree theory,
\begin{equation*}
d(I-\mathcal{T}_{(0,0)}, \Omega, 0)=\pm 1.
\end{equation*}
By (\ref{eq: homotopy invariance}), $d(I-\mathcal{T}_{(1,1)}, \Omega, 0)=\pm 1$. Therefore, the equation $(z,w)=\mathcal{T}_{(1,1)}(z,w)$ has at least one solution, that is, the auxiliary problem with $\lambda=\mu=1$,
\begin{equation*}
\left\{
\begin{split}
& z''(t)-z(t)=f(t,\delta_1(t,z(t)),\delta_2(t,w(t)),z'(t),w'(t)) - \delta_1(t,z(t))\\
& w''(t)-w(t)=g(t,\delta_1(t,z(t)),\delta_2(t,w(t)),z'(t),w'(t)) - \delta_2(t,w(t)),
\end{split}
\right.
\end{equation*}
together with the boundary conditions (\ref{eq: PBC}), has at least one solution $(z_*(t), w_*(t))$.

\bigskip

\textbf{Claim 4.} \textit{The pair $(z_*(t),w_*(t)) \in X^2$, solution of the auxiliary problem (\ref{eq: AP}), (\ref{eq: PBC}), for $\lambda=\mu=1$, is also a solution of the original problem (\ref{eq: OP}), (\ref{eq: PBC}).}

\bigskip

This Claim is proven if the solution $(z_*(t),w_*(t))$ satisfies
\begin{equation}\label{eq: inequalities}
\alpha_1^0(t) \leq z_*(t) \leq \beta_1^0(t), 
\quad
\alpha_2^0(t) \leq w_*(t) \leq \beta_2^0(t),
\quad
\forall t \in [0,T].
\end{equation}

Suppose, by contradiction, that there is $t\in [0,T]$ such that $z_*(t)>\beta_1^0(t)$, and define
\begin{equation}
\max_{t\in[0,T]} \{ z_*(t)-\beta_1^0(t) \}:=z_*(t_0)-\beta_1^0(t_0)>0.
\label{eq: t_0 localization}
\end{equation}

If $t_0 \in ]0,T[$, then  $(z_*-\beta_1^0)'(t_0)=0$ and $(z_*-\beta_1^0)''(t_0)\leq 0$. Then, by (\ref{eq: AP}), (\ref{eq: t_0 localization}), (\ref{eq: f non-increasing}) and (\ref{eq: beta'' < f lower upper}), the following contradiction holds,
\begin{equation}
\begin{split}
0 & \geq z''_*(t_0) - (\beta_1^0)''(t_0) \\
& =
f(t_0, \delta_1(t_0, z_*(t_0)), \delta_2(t_0, w(t_0)), z'_*(t_0), w'(t_0)) - \delta_1(t_0, z_*(t_0)) + z_*(t_0) - \beta_1''(t_0)\\
& =
f(t_0, \beta_1^0(t_0), \delta_2(t_0, w(t_0)), \beta_1'(t_0), w'(t_0)) - \beta_1^0(t_0) + z_*(t_0) - \beta_1''(t_0)\\
& >
f(t_0, \beta_1^0(t_0), \delta_2(t_0, w(t_0)), \beta_1'(t_0), w'(t_0)) -\beta_1''(t_0)\\
& \geq
f(t_0, \beta_1^0(t_0), \beta_2^0(t_0), \beta_1'(t_0), w'(t_0)) -\beta_1''(t_0) \geq 0.
\end{split}
\end{equation}

If $t_0=0$ or $t_0=T$, then $(z_*-\beta_1^0)'(0)\leq 0$ and $(z_*-\beta_1^0)'(T)\geq 0$. By (\ref{eq: beta0}), (\ref{eq: PBC}) and (\ref{eq: beta - monotonicity}),
\begin{equation*}
\begin{split}
0
& \geq (z_*-\beta_1^0)'(0) = (z_*-\beta_1)'(0) = z_*'(T)-\beta_1'(0) \\
& \geq z_*'(T)-\beta_1'(T) = (z_*-\beta_1^0)'(T) \geq 0,
\end{split}
\end{equation*}
so,
\begin{equation}\label{eq: temp-1}
(z_*-\beta_1^0)'(0) = (z_*-\beta_1^0)'(T) = 0,
\end{equation}
and $(z_*-\beta_1^0)''(t_0) \leq 0$. Therefore, we can apply the previous arguments to obtain a similar contradiction, and so, $z_*(t) \leq \beta_1^0(t)$.

Similar arguments can be used to prove the other inequalities in (\ref{eq: inequalities}).
\end{proof}

\bigskip

\begin{example}\label{example}
Consider the following system, for $t\in \lbrack
0,1],$ 
\begin{equation}
\left\{
\begin{split}
& z''(t) = 2 z^3(t) - w(t) + 3 z'(t) - \frac{2}{1+(w'(t))^2} -10  t,\\
& w''(t) = - z(t) + 10 w^3(t) - e^{-(z'(t))^2} - 3 w'(t) - 12 t,
\end{split}
\right.
\label{eq: example - equation}
\end{equation}
together with the periodic boundary conditions (\ref{eq: PBC}).

The functions $\alpha_{i},\beta_{i}:\left[ 0,1\right] \rightarrow \mathbb{R},$ $i=1,2,$ given by
\begin{equation*}
\begin{array}{ccc}
&\alpha_{1}(t) = 1+2t^2-2t^3, \qquad &\beta_{1}(t) = 6/5 - 2t^2 + 2t^3, \\
&\alpha_{2}(t) = 2t - 2t^2, \qquad &\beta_{2}(t) = 1 -3t + 3t^2,
\end{array}
\end{equation*}
\begin{figure}[H]
\begin{subfigure}{.55\textwidth}
  \centering
  \includegraphics[width=1\linewidth]{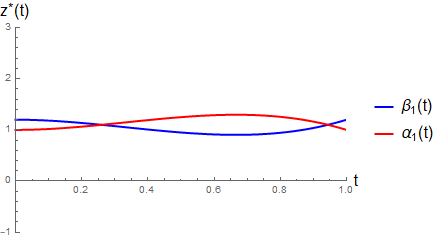}  
  %\caption{$z(t)$}
\end{subfigure}
\begin{subfigure}{.55\textwidth}
  \centering
  \includegraphics[width=1\linewidth]{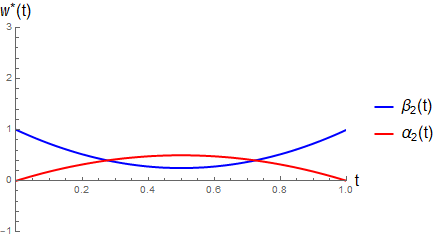}  
  %\caption{$w(t)$}
\end{subfigure}
\caption{Orderless $\alpha_i, \beta_i$ functions, with $i=1,2$.}
\label{fig: example - alpha, beta}
\end{figure}
\noindent are, respectively, lower and upper solutions of problem (\ref{eq: example - equation}), (\ref{eq: PBC}), according to Definition \ref{def: lower upper}, with
\begin{equation*}
\begin{array}{ccc}
&\alpha_{1}^0(t) = -8/27 + 2t^2 - 2t^3, \qquad
 &\beta_{1}^0(t) = 12/5 - 2t^2 + 2t^3, \\
&\alpha_{2}^0(t) = -1/2 +2t - 2t^2, \qquad
 &\beta_{2}^0(t) = 2 -3t + 3t^2.
\end{array}
\end{equation*}

Remark that the lower and upper functions, shown in Figure \ref{fig: example - alpha, beta}, are not ordered, as it is usual in the literature.

The above problem is a particular case of (\ref{eq: OP}), (\ref{eq: PBC}), with $T=1$, and where the non-linearities, according to Definition \ref{def: nagumo}, satisfy a Nagumo-type condition in the set
\begin{equation*}
\widetilde{S} = 
\left\{
\begin{split}
& (t, z_0, w_0, z_1, w_1) \in [0,1]\times \mathbb{R}^4:\\
& \qquad\qquad -8/27 + 2t^2 - 2t^3 \leq z_0 \leq 12/5 - 2t^2 + 2t^3, \\
& \qquad\qquad -1/2 +2t - 2t^2 \leq w_0 \leq 2 -3t + 3t^2 
\end{split}
\right\},
\end{equation*}
such that
\begin{equation*}
\begin{split}
& \left| f(t,z,w,z',w') \right| \leq
2 \left|z\right|^3 + w + 3 \left|z'\right| + \left|\frac{2}{1+(w')^2}\right| + 10\left| t \right|\\
&\qquad\qquad 
\leq 2 \times \frac{12}{5} + 2 + 3 \left|z'\right| + 2 + 10 \\
&\qquad\qquad
= \frac{94}{5} + 3\left|z'\right| := \widetilde{\varphi}(\left| z' \right|),\\
& \left| g(t,z,w,z',w') \right| \leq
\left| z \right| + 10 \left| w \right|^3 + \left| e^ {(z')^2}\right| + 3\left| w' \right| + 12\left| t \right|\\
&\qquad\qquad 
\leq \frac{12}{5} + 10 \times 2^3 + 1 + 3\left| w' \right| + 12\\
&\qquad\qquad 
= \frac{477}{5} + 3\left| w' \right| := \widetilde{\psi}(\left| w' \right|).
\end{split}
\end{equation*}
It is clear that functions $\widetilde{\varphi}$ and $\widetilde{\psi}$ satisfy (\ref{eq: phi,psi}). 

As the assumptions of Theorem \ref{thm: 1} are verified, then the system (\ref{eq: example - equation}), (\ref{eq: PBC}) has, at least, a solution $(z^*(t),w^*(t))\in (C^2[0,1])^{2}$ such that 
\begin{equation}
\begin{split}
-8/27 + 2t^2 - 2t^3 & \leq z^*(t)\leq  12/5 - 2t^2 + 2t^3, \\
-1/2 +2t - 2t^2 & \leq w^*(t)\leq 2 -3t + 3t^2,\quad \forall t\in \lbrack 0,1],
\end{split}
\end{equation}
as shown in Figure \ref{fig: example - alpha0, beta0}.

\begin{figure}[h!]
\begin{subfigure}{.55\textwidth}
  \centering
  \includegraphics[width=1\linewidth]{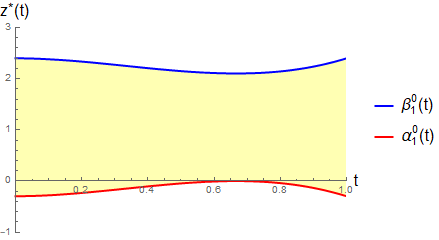}  
  %\caption{$z(t)$}
\end{subfigure}
\begin{subfigure}{.55\textwidth}
  \centering
  \includegraphics[width=1\linewidth]{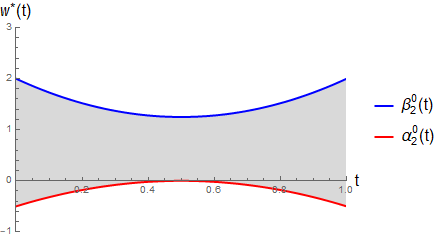}  
  %\caption{$w(t)$}
\end{subfigure}
\caption{Shifted functions, $\alpha_i^0, \beta_i^0, i=1,2$, localizing the solution pair $(z^*(t), w^*(t))$.}
\label{fig: example - alpha0, beta0}
\end{figure}
\end{example}

\bigskip

\section{Coupled forced Van der Pol oscillators}
\label{sec: application}

\quad The equation for the damped harmonic motion is 
\begin{equation*}
x''(t) + \mu x'(t) + x(t) = 0.
\end{equation*}
Balthazar Van der Pol (1889 -- 1959) modified the damped harmonic oscillator by considering a negative quadratic term for the friction term to obtain self-sustained oscillations. This modification resulted in the Van der Pol oscillator \cite{van1926lxxxviii, van1927frequency}, which can be represented by the following equation,
\begin{equation*}
x''(t) - \epsilon (1 - x^2(t)) x'(t) + x(t) = 0,
\end{equation*}
where $x(t)$ is the time-dependent variable and  $-\epsilon (1 - x^2(t))$ is the non-linear damping term, with $\epsilon > 0$.

A variant of this problem can be thought of, by coupling two Van der Pol oscillators. We consider the following system, for $t\in [0,T],$ 
\begin{equation}
\left\{
\begin{split}
& z''(t) = z'(t)(A_1 - B_1 z^2(t)) - C_1 z(t) 
+ D_1 \tanh(E_1 z(t) - F_1 w(t))
+ G_1 \cos(t)\\
& w''(t) = w'(t)(A_2 - B_2 w^2(t)) - C_2 w(t) 
+ D_2 \arctan(E_2 w(t) - F_2 z(t))
+ G_2 \cos(t) 
\end{split}
\right.
\label{eq: app - equation}
\end{equation}
with $A_i, B_i, C_i, D_i, F_i > 0$ and $E_i, G_i \in \mathbb{R}$, $i=1,2$, 
together with the periodic boundary conditions (\ref{eq: PBC}).

The terms $D_1 \tanh(E_1 z(t) - F_1 w(t))$ and $D_2 \arctan(E_2 w(t) - F_2 z(t))$ correspond, respectively, to each coupling and $G_i \cos(t)$ is a time-dependent periodic forcing. We chose $T=1$ and the parameter set
\begin{equation}
\begin{array}{ccccccc}
A_1 = 1 \quad & 
B_1 = 0.5 \quad & 
C_1 = 1 \quad & 
D_1 = 6 \quad & 
E_1 = 3 \quad & 
F_1 = 2 \quad & 
G_1 = 1,\\
A_2 = 1 \quad & 
B_2 = 0.5 \quad & 
C_2 = 1 \quad & 
D_2 = 8 \quad & 
E_2 = 2 \quad & 
F_2 = 1 \quad & 
G_2 = -1.
\end{array}
\label{eq: app - parameter set}
\end{equation}

The functions $\alpha_{i},\beta_{i}:\left[ 0,1\right] \rightarrow \mathbb{R},$ $i=1,2,$ given by
\begin{equation*}
\begin{array}{ccc}
&\alpha_{1}(t) = -1 + t - t^2, \qquad &\beta_{1}(t) = 1 - t^2/2 + t^3/2, \\
&\alpha_{2}(t) = -3/4 + t - t^2, \qquad &\beta_{2}(t) = 1 -t^2 + t^3,
\end{array}
\end{equation*}
are, respectively, lower and upper solutions of problem (\ref{eq: app - equation}), (\ref{eq: PBC}), with the numerical values (\ref{eq: app - parameter set}) according to Definition \ref{def: lower upper}, with
\begin{equation*}
\begin{array}{ccc}
&\alpha_{1}^0(t) = -5/4 + t - t^2, \qquad
 &\beta_{1}^0(t) = 2 -t^2/2 + t^3/2, \\
&\alpha_{2}^0(t) =  -1 + t - t^2, \qquad
 &\beta_{2}^0(t) = 2 -t^2 + t^3.
\end{array}
\end{equation*}

The functions $f$ and $g$ verify the monotonicity requirements and the Nagumo-type condition of Definition \ref{def: nagumo} in the set
\begin{equation*}
S^* =\left\{ 
\begin{split}
& (t, z_0, w_0, z_1, w_1) \in [0,1]\times \mathbb{R}^4:\\
& \qquad\qquad -5/4 + t - t^2 \leq z_0 \leq 2 -t^2/2 + t^3/2,\\
& \qquad\qquad -1 + t - t^2 \leq w_0 \leq 2 -t^2 + t^3 
\end{split}
\right\},
\end{equation*}
with
\begin{equation*}
\begin{split}
& f(t, z_0, w_0, z_1, w_1) \leq
|z_1| \, (1+ 0.5 |z_0|^2) + |z_0| + 6 |\tanh(3z_0 - 2w_0)| + |\cos(t)|\\
& \qquad \qquad
\leq 3|z_1| + 9 := \varphi^*(|z_1|),\\
& g(t, z_0, w_0, z_1, w_1) \leq
|w_1| \, (1+ 0.5 |w_0|^2) + |w_0| + 8 |\arctan(2w_0 - z_0)| + |\cos(t)|\\
& \qquad \qquad
\leq 3|w_1| + 3 + 4\pi := \psi^*(|w_1|).
\end{split}
\end{equation*}

\begin{figure}[h!]
\begin{subfigure}{.55\textwidth}
  \centering
  \includegraphics[width=1\linewidth]{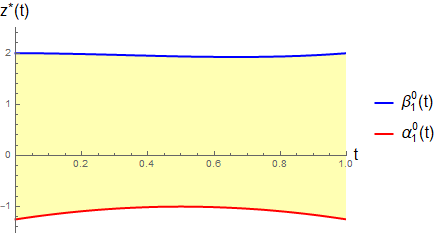}  
  %\caption{$z(t)$}
\end{subfigure}
\begin{subfigure}{.55\textwidth}
  \centering
  \includegraphics[width=1\linewidth]{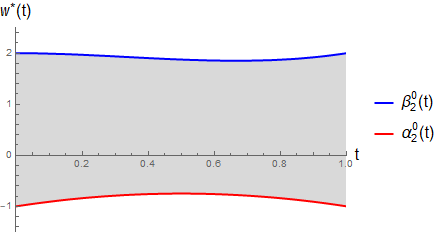}  
  %\caption{$w(t)$}
\end{subfigure}
\caption{Shifted functions, $\alpha_i^0, \beta_i^0, i=1,2$, localizing the solution pair $(z^*(t), w^*(t))$.}
\label{fig: app - alpha0, beta0}
\end{figure}

As all the assumptions of Theorem \ref{thm: 1} are verified, then there is at least a periodic solution for the system of two coupled forced Van der Pol oscilators (\ref{eq: app - equation}) with the periodic boundary conditions (\ref{eq: PBC}), and the parameter set of values (\ref{eq: app - parameter set}). Moreover, remark that, from (\ref{eq: app - equation}) and (\ref{eq: app - parameter set}), this solution is a non-trivial one.

Figure \ref{fig: app - alpha0, beta0} shows the shifted functions that localize existing periodic solutions in the system (\ref{eq: app - equation}), (\ref{eq: PBC}), with parameter values (\ref{eq: app - parameter set}).

\bigskip

\section{Conclusions}
\label{sec: conclusions}

\quad When reviewing the literature, one realizes that there are several attempts to find periodic solutions in second-order differential problems. However, either periodicity requirements are made \cite{lu2019periodic, wang2012ambrosetti}, or approximated solutions are presented \cite{yildirim2011higher, ju2014global, mickens2005generalized}, or the problem does not present complete non-linearities \cite{fonda2021periodic}.

In this work we present a general second-order differential coupled system, with dependencies on both variables, $z,w$, their derivatives, $z',w'$, and on time, $t$. We provide sufficient conditions to prove the existence of periodic solutions for problem (\ref{eq: OP}), (\ref{eq: PBC}) without assuming any type of periodicity in the non-linearities.

To localize an existing periodic solution, no order requirement is made between lower and upper solutions. Therefore, we increase the range of possibilities when considering well-ordered and not well-ordered localizing functions $\alpha_i, \beta_i$, see Example \ref{example}.

Our methodology is successfully applied to a system of coupled Van der Pol oscillators with forcing terms. Under the conditions of the main result, Theorem \ref{thm: 1}, we prove the existence of a periodic solution for problem (\ref{eq: OP}), (\ref{eq: PBC}), and localize it within a strip bounded by shifted lower and upper solutions $\alpha_i^0, \beta_i^0$, whose requirements follow Definition \ref{def: lower upper}.

When applying Theorem \ref{thm: 1} to a specific second-order system, one must verify that the non-linearities satisfy the Nagumo-type condition of Definition \ref{def: nagumo}, as well as the required monotonicities. Note that Theorem \ref{thm: 1} does not guarantee the existence of a \textit{non-trivial} solution. That guarantee arises with the nature of the system in study, \textit{i.e.}, whether if it allows for constant solutions or not.

As we deal with a generalized system, our study has applicability in several real-case scenarios in Nature, and it can be an important tool for other mathematical problems.

\bigskip

\section*{Aknowledgements}

\quad This work research was supported by national funds through the Funda\c{c}\~{a}o para a Ci\^{e}ncia e Tecnologia, FCT, under the project UIDB/04674/2020.

\bigskip

\nocite{*}
\bibliographystyle{unsrt}

\end{document}